\documentclass[11pt]{article}

\usepackage[utf8]{inputenc}
\usepackage[english]{babel}
\usepackage{enumerate}
\usepackage{amsthm,amsmath}
\usepackage{xcolor}
\usepackage{amssymb}
\usepackage{tikz}
\usepackage{aliascnt} 

\usepackage[bookmarks]{hyperref}
\hypersetup{
	pdftitle={Infinitely many $\Delta$-primitive squares},
	pdfauthor={F. Golfieri, V.N. Schvollner},
	colorlinks=true,
	linkcolor=blue,
	citecolor=red,
	filecolor=cyan,
	urlcolor=magenta
}

\usepackage{cleveref}
\newtheorem{theorem}{Theorem}[section]
\newaliascnt{conjecture}{theorem}
\newtheorem{conjecture}[conjecture]{Conjecture}
\aliascntresetthe{conjecture}

\newaliascnt{corollary}{theorem}
\newtheorem{corollary}[corollary]{Corollary}
\aliascntresetthe{corollary}

\newaliascnt{lemma}{theorem}
\newtheorem{lemma}[lemma]{Lemma}
\aliascntresetthe{lemma}

\newaliascnt{problem}{theorem}
\newtheorem{problem}[problem]{Problem}
\aliascntresetthe{problem}

\newaliascnt{proposition}{theorem}
\newtheorem{proposition}[proposition]{Proposition}
\aliascntresetthe{proposition}

\theoremstyle{definition}
\newaliascnt{defini}{theorem}
\newtheorem{defini}[defini]{Definition}
\aliascntresetthe{defini}

\newaliascnt{ejemplo}{theorem}
\newtheorem{ejemplo}[ejemplo]{Example}
\aliascntresetthe{ejemplo}

\newaliascnt{remark}{theorem}
\newtheorem{remark}[remark]{Remark}
\aliascntresetthe{remark}

\crefname{conjecture}{Conjecture}{Conjectures}\Crefname{conjecture}{Conjecture}{Conjectures}
\crefname{problem}{Problem}{Problems}\Crefname{problem}{Problem}{Problems}
\crefname{defini}{Definition}{Definitions}\Crefname{defini}{Definition}{Definitions}
\crefname{ejemplo}{Example}{Examples}\Crefname{ejemplo}{Example}{Examples}
\crefname{observacion}{Observation}{Observations}\Crefname{observacion}{Observation}{Observations}
\title{On the $\Delta$-primitive squares conjecture}
\author{Franco Golfieri\thanks{FAMAF, Universidad Nacional de C\'ordoba (UNC), CONICET, C\'ordoba, Argentina. E-mail: \texttt{francogolfieri87@gmail.com}}
	\and Victor N. Schv\"ollner\thanks{IMASL, CONICET, San Luis, Argentina. E-mail: \texttt{vnsi9m6@gmail.com}}}
\date{}

\begin{document}
	\maketitle
	
	\begin{abstract}
		A natural number $n$ has the $\Delta$ property if some difference of complementary divisors of $n$ is the sum of two such differences; equivalently, by \cite{vnsigma.delta.property}, if some split graph has a factor graph containing an $n$-simple triangle. In additive terms, $n$ has the $\Delta$ property exactly when the positive part of the factor-difference set of $n$, introduced by Erd\H{o}s and Rosenfeld, is not sum-free. The set $\mathbb{N}(\Delta)$ of such numbers is closed under multiplication by squares, and it was conjectured in \cite{vnsigma.delta.property} that it contains infinitely many $\Delta$-primitive perfect squares. We prove this conjecture, and show that these squares can be chosen to avoid any prescribed finite set of primes but $3$.

	\end{abstract}
	
	
	\section{Introduction}
	\label{sec:intro}
	
	For a split graph $S$, the combinatorics of $2$-switches on $S$ is encoded by its factor graph $\Phi(S)$, a multigraph whose induced cycles have length at most $4$ \cite{vnsigma.delta.property}. In \cite{vnsigma.delta.property} the following question was raised: for which $n$ is there a split graph whose factor graph contains an \emph{$n$-simple triangle}, that is, a triangle all of whose edges have multiplicity $n$? The answer is governed by an arithmetic condition on $n$, the $\Delta$ property. Writing $D_n$ for the set of positive divisors of $n$ and
	\[
	D^*_n=\bigl\{\,|a-b|:a,b\in D_n,\ ab=n\,\bigr\},\qquad
	D^+_n=\bigl\{\,x+y:x,y\in D^*_n-\{0\}\,\bigr\},
	\]
	the number $n$ has the \emph{$\Delta$ property}, written $n\in\mathbb{N}(\Delta)$, when $D^*_n\cap D^+_n\neq\varnothing$. Equivalently (\Cref{tool:caracterizacion}), $n$ admits a \emph{$\Delta$-triple}: divisors $1<x<y\leq z<\sqrt{n}$ with
	\[
	\frac{n}{x}-x=\frac{n}{y}-y+\frac{n}{z}-z .
	\]

	\emph{Factor-difference sets.} The set $D^*_n$ has a history of its own: it is the \emph{factor-difference set} of $n$, introduced by Erd\H{o}s and Rosenfeld \cite{er-factor-difference}. Their main question compares the factor-difference sets of different integers: is it true that, for every $s\geq1$, there are integers $L_1<\dots<L_s$ such that $\bigl|\bigcap_{i=1}^{s}D^*_{L_i}\bigr|\geq s$? The case $s=2$ goes back to \cite{er-factor-difference}, and the case $s=3$ is due to Jim\'enez-Urroz \cite{jimenez-urroz}. The problem was revisited by Bremner \cite{bremner}; the answer is known to be positive for $s\leq4$, and the question remains open in general (Problem~885 in \cite{erdosproblems}). Erd\H{o}s and Rosenfeld also studied the small elements of $D^*_n$. If $d\mid n$ and $d\leq\sqrt{n}$, then $h=n/d-d$ satisfies
	\[
	\frac{n}{d}-\sqrt{n}\;\leq\; h\;\leq\; 2\Bigl(\frac{n}{d}-\sqrt{n}\Bigr),
	\]
	because $h=(n/d-\sqrt{n})+(\sqrt{n}-d)$ and $d+n/d\geq2\sqrt{n}$. So the small elements of $D^*_n$ correspond to the divisors of $n$ slightly larger than $\sqrt{n}$. Erd\H{o}s and Rosenfeld proved that infinitely many $n$ have four divisors in $(\sqrt{n},\sqrt{n}+16n^{1/4})$, and that for every $C>0$ all sufficiently large $n$ have at most $1+C^2$ divisors in $[\sqrt{n},\sqrt{n}+Cn^{1/4}]$ \cite{er-factor-difference}; related questions are still open (Problems~886 and~887 in \cite{erdosproblems}). For perfect squares, which are the subject of this paper, Chan \cite{chan2014} proved that every sufficiently large square $n$ has at most five divisors in $\bigl[\sqrt{n}-n^{1/4}(\log n)^{1/7},\,\sqrt{n}+n^{1/4}(\log n)^{1/7}\bigr]$. By the inequality above, if $0<h\leq H$ is an element of $D^*_n$ with $h=n/d-d$ and $d<\sqrt{n}$, then both $d$ and $n/d$ lie in $[\sqrt{n}-H,\sqrt{n}+H]$. Since this interval also contains $\sqrt{n}$ when $n$ is a square, Chan's theorem shows that, for every sufficiently large $N$, the set $D^*_{N^2}$ has at most two elements in $\bigl(0,\,N^{1/2}(2\log N)^{1/7}\bigr]$. Chan also observes that the extremal configurations arise from Pell equations, which is also the source of the squares studied in this paper.

	\emph{The additive viewpoint.} The questions above concern the small elements of a factor-difference set and the way the factor-difference sets of different integers overlap. The $\Delta$ property is of a different nature, since it concerns the additive structure of a single factor-difference set. Indeed, $D^+_n$ is the sumset of $D^*_n-\{0\}$ with itself, and $0\notin D^+_n$. Hence $n\in\mathbb{N}(\Delta)$ if and only if $D^*_n-\{0\}$ is not \emph{sum-free}, that is, if and only if some element of $D^*_n-\{0\}$ is the sum of two, not necessarily distinct, elements of $D^*_n-\{0\}$; see \cite{tao-vu} for a survey on sum-free sets. The arithmetic of this condition is made transparent by the identity $(a+b)^2-(a-b)^2=4ab$, which underlies Fermat's factorization method \cite{pomerance}. If $ab=n$, then $(a-b)^2+4n=(a+b)^2$. Conversely, if $h\geq0$ and $h^2+4n=w^2$ with $w\in\mathbb{N}$, then $w\equiv h\pmod2$, and $(w+h)/2$ and $(w-h)/2$ are complementary divisors of $n$ whose difference is $h$. Therefore
	\begin{equation}
		\label{eq:fermat}
		D^*_n=\bigl\{\,h\in\mathbb{N}\cup\{0\}:\ h^2+4n\ \text{is a perfect square}\,\bigr\},
	\end{equation}
	and $n\in\mathbb{N}(\Delta)$ if and only if there are positive integers $h_1,h_2$ such that $h_1^2+4n$, $h_2^2+4n$ and $(h_1+h_2)^2+4n$ are perfect squares.

	The set $\mathbb{N}(\Delta)$ is closed under multiplication by squares. An element $n\in\mathbb{N}(\Delta)$ is \emph{$\Delta$-primitive} if there are no integers $\alpha,m\geq2$ with $n=\alpha^2m$ and $m\in\mathbb{N}(\Delta)$. In \cite{vnsigma.delta.property} it is shown that there are infinitely many $\Delta$-primitive numbers. The argument uses $F_1(p)=p(2p-1)(3p-2)$ for odd primes $p$: since $p$ divides $F_1(p)$ exactly once, it must divide $m$ in any decomposition $F_1(p)=\alpha^2m$.\footnote{The values $F_1(p)$ themselves need not be $\Delta$-primitive: $F_1(59)=5^2\cdot48321$, and $(39,59,91)$ is a $\Delta$-triple for $48321$.} In a perfect square every prime has even exponent, so this argument is not available, and the following question was left open (\Cref{tool:conj-squares}). In this article, we give a positive answer to that conjecture. Moreover, we prove the following theorem.

	\begin{theorem}
		\label{thm:intro-main}
		There are infinitely many $\Delta$-primitive squares. Moreover, for every finite set $\Pi$ of primes there are infinitely many $\Delta$-primitive squares $\mu^2$ such that no prime factor of $\mu$ other than $3$ lies in $\Pi$.
	\end{theorem}

	The squares come from the $\Delta$-triple polynomial family
	\[
	\mathcal{F}_1:\ \bigl(t,\ 2t-1,\ 2t+1\bigr),
	\]
	whose square values $N_k^2$ are governed by the Pell equation $u^2-3v^2=1$ (\Cref{thm:F1}). This is possible because $n_1(t)=t^2\cdot3(4t^2-1)$ has the repeated factor $t^2$, so that $n_1(t)$ is a square exactly when $3(4t^2-1)$ is, which is a Pell condition. By contrast, if $f=R^2S$ with $S$ squarefree of degree at least $3$, then Siegel's theorem on integral points \cite{siegel1929} allows only finitely many square values of $f$, an obstruction already discussed in \cite[Section~4.4]{vnsigma.delta.property}. The triples of $\mathcal{F}_1$ are \emph{generic}, i.e.\ $y<z$, and this is unavoidable: we show that a perfect square has no \emph{duplicated} $\Delta$-triple, i.e.\ none with $y=z$ (\Cref{thm:no-dup-squares}), by reducing the question to a theorem of Euler on products of four terms in arithmetic progression.

	\emph{Squares and Pythagorean triangles.} For perfect squares, the main subject of this paper, \eqref{eq:fermat} has a classical geometric meaning. Note first that $0\in D^*_n$ if and only if $n$ is a square; if $0$ were allowed in the definition of $D^+_n$, every square would trivially have the $\Delta$ property. For $n=N^2$, \eqref{eq:fermat} says that the positive elements of $D^*_{N^2}$ are exactly the legs $h$ of the Pythagorean triangles whose other leg is $2N$, that is, the positive integers $h$ such that $h^2+(2N)^2$ is a perfect square. Consequently, $N^2\in\mathbb{N}(\Delta)$ if and only if there are three Pythagorean triangles with the common leg $2N$ whose other legs $h_1,h_2,h_3$ satisfy $h_1+h_2=h_3$. For instance, the smallest $\Delta$-primitive square is $30^2$ (\Cref{ex:min-S}), and its only $\Delta$-triple $(9,10,25)$ corresponds to the triangles
	\[
	(11,60,61),\qquad (80,60,100),\qquad (91,60,109),
	\]
	where $11+80=91$. Geometrically, the $91\times60$ rectangle has the integral diagonal $109$, and the point of one of its sides of length $91$ at distance $11$ from a vertex lies at the integral distances $61$ and $100$ from the two vertices of the opposite side. A configuration with common leg $2d$ can be scaled to one with common leg $2N$ whenever $d\mid N$. By \Cref{prop:min-S}, $N^2$ is $\Delta$-primitive if and only if the leg $2N$ carries such a configuration, but no leg $2d$ with $d$ a proper divisor of $N$ does. In this language, \Cref{thm:intro-main} provides infinitely many such $N$.


	
	
	\section{Preliminaries}
	\label{sec:tools}
	
	We use the notation of \cite{vnsigma.delta.property}, and write $\mathbb{N}(\Delta^2)$ for the set of perfect squares in $\mathbb{N}(\Delta)$.
	
	\begin{proposition}[\cite{vnsigma.delta.property}]
		\label{tool:caracterizacion}
		A natural number $n$ lies in $\mathbb{N}(\Delta)$ if and only if there exists $(x,y,z)\in D_n^3$ with $1<x<y\leq z<\sqrt{n}$ and
		\begin{equation}
			\frac{n}{x}-x \;=\; \frac{n}{y}-y + \frac{n}{z}-z.
			\label{eq:tool-delta}
		\end{equation}
		Such a triple is a \emph{$\Delta$-triple} for $n$. Equation \eqref{eq:tool-delta} is equivalent to
		\begin{equation}
			(xy+xz-yz)\,n \;=\; xyz\,(z+y-x).
			\label{eq:tool-delta2}
		\end{equation}
	\end{proposition}
	
	\begin{proposition}[\cite{vnsigma.delta.property}]
		\label{tool:rescaling}
		If $n\in\mathbb{N}(\Delta)$, then $\alpha^2 n\in\mathbb{N}(\Delta)$ for all $\alpha\in\mathbb{N}$. Every $n\in\mathbb{N}(\Delta)$ is either $\Delta$-primitive or of the form $\alpha^2m$ with $\alpha\geq2$ and $m$ $\Delta$-primitive. Moreover $(x,2x-1,2x-1)$ is a $\Delta$-triple for $F_1(x)=x(2x-1)(3x-2)$ for every $x\geq2$.
	\end{proposition}
	
	\begin{conjecture}[\cite{vnsigma.delta.property}]
		\label{tool:conj-squares}
		There are infinitely many $\Delta$-primitive squares.
	\end{conjecture}
	
    Let $u_k+v_k\sqrt3=(2+\sqrt3)^k$ for $k\geq0$. The pairs $(u_k,v_k)$, $k\geq1$, are the positive solutions of
	\begin{equation}
		\label{eq:pell}
		u^2-3v^2=1,
	\end{equation}
	and $u_k$ is even if and only if $k$ is odd. The following facts are classical \cite{lucas1878}:
	\begin{itemize}
		\item $v_{2k}=2u_kv_k$, as the coefficient of $\sqrt3$ in $(u_k+v_k\sqrt3)^2$;
		\item the sequence $(v_k)$ has the strong divisibility property
		\begin{equation}
			\label{eq:strong}
			\gcd(v_j,v_k)=v_{\gcd(j,k)};
		\end{equation}
		\item for a prime $\ell\geq5$, the least $r\geq1$ with $\ell\mid v_r$ exists and divides $\ell-\bigl(\frac3\ell\bigr)$; moreover $\ell\mid v_m$ if and only if $r\mid m$.
	\end{itemize}

A $\Delta$-triple $(x,y,z)$ is \emph{duplicated} if $y=z$, and \emph{generic} if $y<z$. We need two more facts.

\begin{proposition}[{\cite[Proposition~2.4]{vnsigma.delta.property}}]
	\label{tool:y/x}
	If $(x,y,z)$ is a $\Delta$-triple for $n$, then $1<y/x<2$.
\end{proposition}

\begin{theorem}[Euler; see {\cite[pp.~440 and 635]{dickson}} and \cite{hlst}]
	\label{tool:euler}
	If $c,\delta$ are coprime positive integers, then $c(c+\delta)(c+2\delta)(c+3\delta)$ is not a perfect square.
\end{theorem}

The following result shows that the duplicated regime plays no role for squares.

\begin{theorem}
	\label{thm:no-dup-squares}
	A perfect square has no duplicated $\Delta$-triple.
\end{theorem}

\begin{proof}
	Let $(x,y,y)$ be a duplicated $\Delta$-triple for $n$. Putting $z=y$ in \eqref{eq:tool-delta2} gives $(2x-y)\,n=xy\,(2y-x)$. By \Cref{tool:y/x}, $x<y<2x$, so $c=2x-y$ and $\delta=y-x$ are positive integers, and $x=c+\delta$, $y=c+2\delta$, $2y-x=c+3\delta$. Multiplying by $c$, we get
	\[
	c^2n=c(c+\delta)(c+2\delta)(c+3\delta).
	\]
	Let $g=\gcd(c,\delta)$, and write $c=gc'$ and $\delta=g\delta'$. Then
	\[
	c'(c'+\delta')(c'+2\delta')(c'+3\delta')=\frac{c^2n}{g^4}.
	\]
	If $n$ were a perfect square, the right-hand side would be the square of a rational number, so the integer on the left-hand side would be a perfect square. Since $\gcd(c',\delta')=1$, this contradicts \Cref{tool:euler}.
\end{proof}

So all the $\Delta$-triples considered below are generic.

In terms of \eqref{eq:fermat}, a duplicated $\Delta$-triple $(x,y,y)$ for $n$ yields a positive integer $h$ such that $h$ and $2h$ lie in $D^*_n$, namely $h=n/y-y$ and $2h=n/x-x$. Conversely, if $h\geq1$ and $h,2h\in D^*_n$, then $2h=h+h\in D^*_n\cap D^+_n$, and the proof of \Cref{tool:caracterizacion} in \cite{vnsigma.delta.property} provides a $\Delta$-triple $(x,y,z)$ with $n/y-y=n/z-z=h$; hence $y=z$, because $d\mapsto n/d-d$ is strictly decreasing on $(0,\sqrt{n}\,]$. Since $(2h)^2+4n=4(h^2+n)$, it follows that $n$ has a duplicated $\Delta$-triple if and only if $h^2+n$ and $h^2+4n$ are perfect squares for some $h\geq1$. For $n=N^2$, \Cref{thm:no-dup-squares} thus says that no positive integer is a leg of two Pythagorean triangles whose other legs are $N$ and $2N$, respectively.

Our squares come from the following family.

\begin{theorem}
	\label{thm:F1}
	For every integer $t\geq2$, $(t,\,2t-1,\,2t+1)$ is a $\Delta$-triple for $n_1(t)=3t^2(4t^2-1)$. Moreover, $n_1(t)$ is a perfect square if and only if $2t=u_k$ for some odd $k$, and then $n_1(t)=(3tv_k)^2$.
\end{theorem}

\begin{proof}
	The quotients $n_1(t)/t=3t(2t-1)(2t+1)$, $n_1(t)/(2t-1)=3t^2(2t+1)$ and $n_1(t)/(2t+1)=3t^2(2t-1)$ are integers. The inequalities $1<t<2t-1<2t+1<\sqrt{n_1(t)}$ hold for $t\geq2$, since $n_1(t)-(2t+1)^2=12t^4-7t^2-4t-1>0$. Finally,
	\begin{align*}
		\frac{n_1(t)}{t}-t=12t^3-4t&=\bigl(6t^3+3t^2-2t+1\bigr)+\bigl(6t^3-3t^2-2t-1\bigr)\\
		&=\Bigl(\frac{n_1(t)}{2t-1}-(2t-1)\Bigr)+\Bigl(\frac{n_1(t)}{2t+1}-(2t+1)\Bigr),
	\end{align*}
	so \eqref{eq:tool-delta} holds, and $(t,2t-1,2t+1)$ is a $\Delta$-triple for $n_1(t)$ by \Cref{tool:caracterizacion}.

	Since $n_1(t)=t^2\cdot3(4t^2-1)$, comparing $\ell$-adic valuations shows that $n_1(t)$ is a perfect square if and only if $3(4t^2-1)=w^2$ for some integer $w$. Then $3\mid w$, and writing $w=3v$ gives $(2t)^2-3v^2=1$; conversely, every solution $(2t,v)$ of \eqref{eq:pell} gives $3(4t^2-1)=(3v)^2$. Since $2t$ is even, the positive solutions of \eqref{eq:pell} of this form are exactly the pairs $(u_k,v_k)$ with $k$ odd, and then $n_1(t)=t^2(3v_k)^2=(3tv_k)^2$.
\end{proof}

\section{Infinitely many $\Delta$-primitive squares}
	\label{sec:proof}
	
	We recast the conjecture as a statement about divisibility among square roots, derive a sufficient criterion, and verify it with the square values of $\mathcal{F}_1$.
	
	\begin{defini}
		\label{def:S}
		Let
		\[
		\mathcal{S}=\{\,N\in\mathbb{N}:\ N^2\in\mathbb{N}(\Delta)\,\},
		\]
		so that $N\mapsto N^2$ maps $\mathcal{S}$ bijectively onto $\mathbb{N}(\Delta^2)$. An element $N\in\mathcal{S}$ is \emph{minimal} if no divisor $d$ of $N$ with $d<N$ belongs to $\mathcal{S}$. We write $\operatorname{Min}(\mathcal{S})$ for the set of minimal elements of $\mathcal{S}$.
	\end{defini}
	
	The first basic fact is that $\mathcal{S}$ is closed under taking multiples.
	
	\begin{lemma}
		\label{lem:S-up}
		If $N\in\mathcal{S}$ and $N\mid M$, then $M\in\mathcal{S}$.
	\end{lemma}
	
	\begin{proof}
		Write $M=\alpha N$ with $\alpha\in\mathbb{N}$. Then $M^2=\alpha^2N^2\in\mathbb{N}(\Delta)$ by \Cref{tool:rescaling}.
	\end{proof}
	
	The second basic fact identifies the $\Delta$-primitive squares inside $\mathcal{S}$.
	
	\begin{proposition}
		\label{prop:min-S}
		Let $N\in\mathbb{N}$. Then $N^2$ is $\Delta$-primitive if and only if $N\in\operatorname{Min}(\mathcal{S})$. Consequently, $N\mapsto N^2$ is a bijection from $\operatorname{Min}(\mathcal{S})$ onto the set of $\Delta$-primitive squares, and \Cref{tool:conj-squares} is equivalent to the statement that $\operatorname{Min}(\mathcal{S})$ is infinite.
	\end{proposition}
	
	\begin{proof}
		Note first that $1\notin\mathbb{N}(\Delta)$. Suppose that $N^2$ is $\Delta$-primitive, so $N\in\mathcal{S}$. Let $d$ be a divisor of $N$ with $d<N$, and suppose that $d\in\mathcal{S}$. Put $\alpha=N/d$ and $m=d^2$. Then $\alpha\geq 2$, $m\in\mathbb{N}(\Delta)$, $m\geq 2$ (as $1\notin\mathbb{N}(\Delta)$), and $N^2=\alpha^2m$, which contradicts the $\Delta$-primitivity of $N^2$. Hence no such $d$ exists, that is, $N\in\operatorname{Min}(\mathcal{S})$.
		
		Conversely, let $N\in\operatorname{Min}(\mathcal{S})$. Then $N^2\in\mathbb{N}(\Delta)$. Suppose that $N^2=\alpha^2m$ for some integers $\alpha,m\geq 2$ with $m\in\mathbb{N}(\Delta)$. Since $\alpha^2\mid N^2$, comparing $\ell$-adic valuations gives $\alpha\mid N$. Put $d=N/\alpha$. Then $d$ is a divisor of $N$, $d<N$ because $\alpha\geq 2$, and $d^2=N^2/\alpha^2=m\in\mathbb{N}(\Delta)$, so $d\in\mathcal{S}$. This contradicts the minimality of $N$. Hence $N^2$ is $\Delta$-primitive.
		
		The remaining assertions follow, since $N\mapsto N^2$ is injective on $\mathbb{N}$.
	\end{proof}
	
	
	\begin{corollary}
		\label{cor:min-test}
		Let $N\in\mathcal{S}$. Then $N\in\operatorname{Min}(\mathcal{S})$ if and only if $N/q\notin\mathcal{S}$ for every prime divisor $q$ of $N$.
	\end{corollary}
	
	\begin{proof}
		If $N$ is minimal, then $N/q\notin\mathcal{S}$ for every prime $q\mid N$, since $N/q$ is a divisor of $N$ smaller than $N$. Conversely, assume $N/q\notin\mathcal{S}$ for every prime $q\mid N$, and let $d$ be a divisor of $N$ with $d<N$. Then $N/d>1$; let $q$ be a prime divisor of $N/d$, which is also a prime divisor of $N$. Since $N/q=d\cdot\frac{N/d}{q}$, we have $d\mid N/q$. If $d$ belonged to $\mathcal{S}$, then $N/q$ would belong to $\mathcal{S}$ by \Cref{lem:S-up}, a contradiction. Hence $d\notin\mathcal{S}$, and $N$ is minimal.
	\end{proof}
	
	\begin{ejemplo}
		\label{ex:min-S}
		A computer search shows that the elements of $\operatorname{Min}(\mathcal{S})$ below $1500$ are
		\[
		30,\ 84,\ 140,\ 234,\ 315,\ 440,\ 455,\ 462,\ 585,\ 715,\ 765,\ 819,\ 1155 ;
		\]
		by \Cref{prop:min-S}, their squares are exactly the $\Delta$-primitive squares below $1500^2$. In particular, the two smallest $\Delta$-primitive squares are $30^2$ and $84^2$, in agreement with \cite{vnsigma.delta.property}; and $585^2=n_1(13)$ is the smallest square value $n_1(t_0)$ with $t_0\geq 2$.
		
	\end{ejemplo}

	\begin{lemma}
		\label{lem:below-min}
		For every $N\in\mathcal{S}$ there exists $\mu\in\operatorname{Min}(\mathcal{S})$ with $\mu\mid N$.
	\end{lemma}
	
	\begin{proof}
		The set of divisors of $N$ that belong to $\mathcal{S}$ is finite and contains $N$; let $\mu$ be its least element. If $d$ is a divisor of $\mu$ with $d<\mu$, then $d$ is also a divisor of $N$, and therefore $d\notin\mathcal{S}$ by the choice of $\mu$. Hence $\mu\in\operatorname{Min}(\mathcal{S})$.
	\end{proof}
	
	\begin{theorem}
		\label{thm:criterion}
		Let $A\subseteq\mathcal{S}$ be a set such that $\gcd(a,b)\notin\mathcal{S}$ for all distinct $a,b\in A$. Then $\operatorname{Min}(\mathcal{S})$ has at least $|A|$ elements. In particular, if $A$ is infinite, then there are infinitely many $\Delta$-primitive squares.
	\end{theorem}
	
	\begin{proof}
		By \Cref{lem:below-min}, we may choose for each $a\in A$ an element $\mu(a)\in\operatorname{Min}(\mathcal{S})$ with $\mu(a)\mid a$. Let $a,b\in A$ be distinct, and suppose that $\mu(a)=\mu(b)$. Then $\mu(a)$ divides both $a$ and $b$, hence divides $\gcd(a,b)$, and \Cref{lem:S-up} gives $\gcd(a,b)\in\mathcal{S}$, contrary to the hypothesis. Therefore $a\mapsto\mu(a)$ is an injective map from $A$ to $\operatorname{Min}(\mathcal{S})$, which proves the first assertion. The second one follows from \Cref{prop:min-S}.
	\end{proof}
	
	
	
	Let $(u_k,v_k)_{k\geq 0}$ be as in \Cref{sec:tools}. For odd $k$, the integer $u_k$ is even, and we put
	\[
	t_k=\frac{u_k}{2}, \qquad N_k=3\,t_k\,v_k \qquad (k\ \text{odd}).
	\]
	The first values are $N_1=3$, $N_3=585$, $N_5=113487$ and $N_7=22015893$.
	
	\begin{lemma}
		\label{lem:Nk}
		Let $k\geq 1$ be odd. Then:
		\begin{enumerate}
			\item $N_k^2=n_1(t_k)$, where $n_1=3t^2(4t^2-1)$ is the polynomial of \Cref{thm:F1};
			\item $4N_k=3\,v_{2k}$;
			\item if $k\geq 3$, then $N_k\in\mathcal{S}$.
		\end{enumerate}
	\end{lemma}
	
	\begin{proof}
		(1) Since $u_k=2t_k$ and $u_k^2-3v_k^2=1$, we have $4t_k^2-1=3v_k^2$. Hence
		\[
		n_1(t_k)=3t_k^2(4t_k^2-1)=9t_k^2v_k^2=N_k^2 .
		\]
		
		(2) We have $v_{2k}=2u_kv_k=4t_kv_k$, so $3v_{2k}=12t_kv_k=4N_k$.
		
		(3) The sequence $(u_k)_{k\geq 1}$ is strictly increasing and $u_3=26$, so $t_k\geq 13$ for every odd $k\geq 3$. By \Cref{thm:F1}, applied with $t=t_k\geq 2$, the triple $(t_k,2t_k-1,2t_k+1)$ is a $\Delta$-triple for $n_1(t_k)$, so $n_1(t_k)\in\mathbb{N}(\Delta)$ by \Cref{tool:caracterizacion}. By (1), $N_k^2\in\mathbb{N}(\Delta)$, that is, $N_k\in\mathcal{S}$.
	\end{proof}
	
	\begin{lemma}
		\label{lem:gcd-Nk}
		For all odd $m,k\geq 1$,
		\[
		\gcd(N_m,N_k)=N_{\gcd(m,k)} .
		\]
		In particular, $\gcd(N_m,N_k)=N_1=3$ whenever $\gcd(m,k)=1$, and $N_a\mid N_k$ whenever $a\mid k$.
	\end{lemma}
	
	\begin{proof}
		Put $d=\gcd(m,k)$, so that $\gcd(2m,2k)=2d$. Since $d$ is odd, $N_d$ is defined. By \Cref{lem:Nk}(2), applied to each of the indices $m$, $k$ and $d$, and by \eqref{eq:strong}, we obtain
		\begin{align*}
			4\gcd(N_m,N_k)&=\gcd(4N_m,4N_k)=\gcd(3v_{2m},3v_{2k})\\
			&=3\gcd(v_{2m},v_{2k})=3v_{2d}=4N_d,
		\end{align*}
		proving the lemma.
	\end{proof}
	
	%
	
	\begin{theorem}
		\label{thm:main}
		There are infinitely many $\Delta$-primitive squares. That is, \Cref{tool:conj-squares} holds.
	\end{theorem}
	
	
	\begin{proof}
		Let $A=\{N_p:\ p\ \text{an odd prime}\}$. By \Cref{lem:Nk}(3), $A\subseteq\mathcal{S}$. Since $(u_k)_{k\geq 1}$ and $(v_k)_{k\geq 1}$ are strictly increasing sequences of positive integers, so is $(N_k)_{k\ \text{odd}}$; hence $A$ is infinite. Let $p\neq p'$ be odd primes. Then $\gcd(p,p')=1$, so $\gcd(N_p,N_{p'})=N_1=3$ by \Cref{lem:gcd-Nk}. Moreover $3\notin\mathcal{S}$: the divisors of $9$ are $1$, $3$ and $9$, none of which lies strictly between $1$ and $\sqrt{9}=3$, so $9$ has no $\Delta$-triple and $9\notin\mathbb{N}(\Delta)$ by \Cref{tool:caracterizacion}. Thus $\gcd(a,b)\notin\mathcal{S}$ for all distinct $a,b\in A$, and \Cref{thm:criterion} shows that there are infinitely many $\Delta$-primitive squares.
	\end{proof}
	
	In words: each $N_p$ is the square root of an element of $\mathbb{N}(\Delta^2)$, and therefore lies above the square root $\mu_p$ of some $\Delta$-primitive square; two distinct $N_p$ have only the factor $3$ in common, and $3$ is too small to be the square root of an element of $\mathbb{N}(\Delta)$; hence the $\mu_p$ are pairwise distinct. Note that any infinite set of pairwise coprime odd integers $k\geq 3$ could replace the odd primes in the proof.
	
	\begin{theorem}
	\label{thm:avoid}
	Let $p\geq5$ be prime, and let $\mu_p\in\operatorname{Min}(\mathcal{S})$ divide $N_p$. Every prime factor of $\mu_p$ other than $3$ is $\equiv\pm1\pmod p$, and in particular is at least $2p-1$. Consequently, for every finite set $\Pi$ of primes there are infinitely many $\Delta$-primitive squares $\mu^2$ such that no prime factor of $\mu$ other than $3$ lies in $\Pi$.
\end{theorem}


\begin{proof}
	We first show that $N_p=3t_pv_p$ is odd. Since $u_p=2t_p$ and $u_p^2-3v_p^2=1$, we have $3v_p^2=4t_p^2-1$, which is odd; hence $v_p$ is odd. Then $v_p^2\equiv1\pmod 8$, so $4t_p^2=3v_p^2+1\equiv4\pmod 8$, and $t_p$ is odd as well. Next, $3$ divides $N_p$ exactly once. Indeed, $u_p\equiv2^p\equiv2\pmod 3$, so $t_p\equiv1\pmod3$; and $3\nmid v_p$ by \eqref{eq:strong}, since $\gcd(v_p,v_3)=v_{\gcd(p,3)}=v_1=1$ and $3\mid v_3=15$.
	
	Now let $\ell\neq3$ be a prime dividing $N_p$. Since $N_p$ is odd, $\ell\geq5$. As $4N_p=3v_{2p}$, we have $\ell\mid v_{2p}$, so the least $r\geq1$ with $\ell\mid v_r$ divides $2p$. Moreover $r\notin\{1,2\}$, because $v_1=1$ and $v_2=4$. Hence $r\in\{p,2p\}$, and since $r$ divides $\ell-\bigl(\frac3\ell\bigr)$, we get $p\mid\ell-\bigl(\frac3\ell\bigr)$, that is, $\ell\equiv\pm1\pmod p$. The numbers $p-1$ and $p+1$ are even, so $\ell\notin\{p-1,p+1\}$, and therefore $\ell\geq2p-1$. Since $\mu_p\mid N_p$, every prime factor of $\mu_p$ other than $3$ has these properties, which proves the first assertion.
	
   For the last assertion, take the primes $p$ with $2p-1>\max\Pi$; the corresponding $\mu_p$ are pairwise distinct by the proof of \Cref{thm:main}.
\end{proof}

\begin{remark}
		\label{obs:not-identified}
		The proof produces, for every odd prime $p$, a $\Delta$-primitive square $\mu_p^2$ with $\mu_p\mid N_p$, the $\mu_p$ being pairwise distinct, but does not identify $\mu_p$. A computer verification based on \Cref{cor:min-test} shows that $\mu_p=N_p$, i.e.\ $N_p^2$ is itself $\Delta$-primitive, for every odd prime $p\leq71$ and for $p\in\{79,89,97,101,103,107,109,137,139\}$. Composite indices behave differently: if $k$ is odd with a divisor $a$, $3\leq a<k$, then $N_a\in\mathcal{S}$, $N_a\mid N_k$ and $N_a<N_k$, so $N_k\notin\operatorname{Min}(\mathcal{S})$. For instance $N_3=585$ divides $N_9=4\,270\,969\,755$.
	\end{remark}
	
	\begin{problem}
		\label{prob:prime-index}
		Is $N_p^2$ $\Delta$-primitive for every odd prime $p$? For infinitely many odd primes $p$?
	\end{problem}

A possible approach to study \Cref{prob:prime-index} is by looking at all $\Delta$-triples of $N_k^2$, not only at the one coming from $\mathcal{F}_1$.

\begin{lemma}
	\label{lem:rho-sigma}
	For odd $k$ define $\rho_k,\sigma_k$ by $\rho_k+\sigma_k\sqrt3=(1+\sqrt3)(2+\sqrt3)^{(k-1)/2}$. Then:
	\begin{enumerate}
		\item the pairs $(\rho_k,\sigma_k)$ are the positive solutions of $\rho^2-3\sigma^2=-2$;
		\item $2t_k-1=\rho_k^2$, $2t_k+1=3\sigma_k^2$ and $v_k=\rho_k\sigma_k$, so that $N_k=3\,t_k\,\rho_k\,\sigma_k$ and the $\Delta$-triple of $\mathcal{F}_1$ for $N_k^2$ is $(t_k,\rho_k^2,3\sigma_k^2)$;
		\item $t_k$, $\rho_k$ and $\sigma_k$ are pairwise coprime, $3\nmid t_k\rho_k$, and $3\nmid\sigma_k$ if $3\nmid k$.
	\end{enumerate}
\end{lemma}

\begin{proof}
	\emph{(1)} Since $1+\sqrt3$ has norm $-2$, the pairs $(\rho_k,\sigma_k)$ solve $\rho^2-3\sigma^2=-2$, and the classical descent for \eqref{eq:pell} shows that there are no other positive solutions.

	\emph{(2)} Since $(1+\sqrt3)^2=2(2+\sqrt3)$, we have $(\rho_k+\sigma_k\sqrt3)^2=2u_k+2v_k\sqrt3$. Hence $\rho_k^2+3\sigma_k^2=4t_k$ and $\rho_k\sigma_k=v_k$, and (2) follows using $\rho_k^2-3\sigma_k^2=-2$.

	\emph{(3)} Coprimality follows from $2t_k-\rho_k^2=1$, $3\sigma_k^2-2t_k=1$ and the fact that $\rho_k$ is odd. Modulo $3$ we have $2t_k\equiv-1$ and $\rho_k^2\equiv1$. If $3\nmid k$, then \eqref{eq:strong} gives $\gcd(v_k,v_3)=v_1=1$, and since $v_3=15$, $3\nmid v_k=\rho_k\sigma_k$.
\end{proof}

\begin{proposition}
	\label{prop:common-factor}
	Let $N\in\mathcal{S}$. If no $\Delta$-triple for $N^2$ has its three components divisible by a common prime factor of $N$, then $N^2$ is $\Delta$-primitive.
\end{proposition}

\begin{proof}
	Suppose that $N^2$ is not $\Delta$-primitive. By \Cref{prop:min-S} and \Cref{cor:min-test}, there is a prime $q\mid N$ with $N/q\in\mathcal{S}$. Let $(x',y',z')$ be a $\Delta$-triple for $(N/q)^2$. Since $x'\mid(N/q)^2$, we have $qx'\mid q(N/q)^2$, which divides $N^2$; similarly for $y'$ and $z'$. Moreover $1<qx'<qy'\leq qz'<q\cdot N/q=N$, and
	\[
	\frac{N^2}{qx'}-qx'=q\Bigl(\frac{(N/q)^2}{x'}-x'\Bigr),
	\]
	and likewise for $y'$ and $z'$. Hence \eqref{eq:tool-delta} for $(x',y',z')$ and $(N/q)^2$ gives \eqref{eq:tool-delta} for $(qx',qy',qz')$ and $N^2$, so $(qx',qy',qz')$ is a $\Delta$-triple for $N^2$ whose components are all divisible by the prime factor $q$ of $N$.
\end{proof}

\begin{corollary}
	\label{cor:unique-primitive}
	Let $k\geq3$ be odd. If every $\Delta$-triple for $N_k^2$ has coprime components, then $N_k^2$ is $\Delta$-primitive. In particular, this holds if $(t_k,\rho_k^2,3\sigma_k^2)$ is the only $\Delta$-triple for $N_k^2$.
\end{corollary}

\begin{proof}
	The first assertion is a special case of \Cref{prop:common-factor}. For the second, the components of $(t_k,\rho_k^2,3\sigma_k^2)$ are coprime, since $\gcd(t_k,\rho_k)=1$ by \Cref{lem:rho-sigma}(3).
\end{proof}

\begin{remark}
	\label{obs:evidence}
	A complete enumeration shows that $(t_p,\rho_p^2,3\sigma_p^2)$ is the only $\Delta$-triple for $N_p^2$ for every odd prime $p\leq71$, and also for $p\in\{79,89,97,101,103,107,109,137,139\}$. By \Cref{cor:unique-primitive}, all these $N_p^2$ are $\Delta$-primitive.

	Composite indices behave differently. Besides the multiples of the triples of $\mathcal{F}_1$ for $N_d^2$ with $d\mid k$, the square $N_9^2$ has two further $\Delta$-triples, coming from $765^2$ and from $104463^2$, and $N_{15}^2$ has eight.
\end{remark}

\begin{conjecture}
	\label{conj:unique}
	For every odd prime $p$, the only $\Delta$-triple for $N_p^2$ is $(t_p,\rho_p^2,3\sigma_p^2)$.
\end{conjecture}

By \Cref{cor:unique-primitive}, \Cref{conj:unique} implies an affirmative answer to \Cref{prob:prime-index}. 


\section{$\Delta$-primitive numbers whose square has the $\Delta$ property}
\label{sec:n-and-square}

If $n\in\mathbb{N}(\Delta)$ and $n^2\in\mathbb{N}(\Delta)$, then $n^k\in\mathbb{N}(\Delta)$ for every $k\geq1$, by \Cref{tool:rescaling}, since $n^{2j+1}=(n^j)^2\,n$ and $n^{2j+2}=(n^j)^2\,n^2$; this was already observed in \cite[Section~3]{vnsigma.delta.property}. Hence $n$ and $n^2$ have the $\Delta$ property if and only if every power of $n$ does, and this section is about the numbers all of whose powers lie in $\mathbb{N}(\Delta)$. By \Cref{tool:rescaling} again, the set of such $n$ is closed under multiplication by squares, so it is infinite. Its least element is $60$, with the $\Delta$-triples $(2,3,4)$ for $60$ and $(18,20,50)$ for $60^2$. The question becomes interesting once primitivity is required. We shall use the following classical result.

\begin{theorem}[Dirichlet \cite{dirichlet1837}]
	\label{tool:dirichlet}
	If $a$ and $q$ are coprime positive integers, then there are infinitely many primes $p\equiv a\pmod q$.
\end{theorem}

\begin{theorem}
	\label{thm:n-and-square}
	There are infinitely many $\Delta$-primitive numbers $n$ such that $n^2\in\mathbb{N}(\Delta)$.
\end{theorem}

\begin{proof}
   Since $228=2^2\cdot3\cdot19$ and $455=5\cdot7\cdot13$ are coprime, \Cref{tool:dirichlet} gives infinitely many primes $p\equiv228\pmod{455^2}$. For such $p$ we have $2p-1=455\,(1+910j)$ for some $j$, and $910=2\cdot5\cdot7\cdot13$. Hence $5$, $7$ and $13$ divide $2p-1$ exactly once. They do not divide $3p-2\equiv682\pmod{455}$, and $p\nmid(2p-1)(3p-2)$. So each of $5,7,13,p$ divides $F_1(p)$ exactly once.
 By \Cref{tool:rescaling}, $F_1(p)=\alpha^2m$ with $\alpha\geq1$ and $m$ $\Delta$-primitive. Each $\ell\in\{5,7,13,p\}$ divides $F_1(p)$ exactly once, so $\ell\nmid\alpha$ and hence $\ell\mid m$. In particular $455\mid m$. Since $455\in\operatorname{Min}(\mathcal{S})$ (\Cref{ex:min-S}), \Cref{lem:S-up} gives $m^2\in\mathbb{N}(\Delta)$. Finally $p\leq m\leq F_1(p)$. Choosing primes with $p_{i+1}>F_1(p_i)$ gives strictly increasing values of $m$.
\end{proof}

For instance, $p=683$ gives the squarefree number $F_1(683)=3\cdot5\cdot7\cdot13\cdot23\cdot89\cdot683$, which is $\Delta$-primitive and whose square lies in $\mathbb{N}(\Delta)$.

\begin{remark}
	\label{rem:why-455}
	The numbers $455$ and $228$ in the proof of \Cref{thm:n-and-square} are the least choices compatible with the method. The proof needs some $\mu\in\mathcal{S}$ with $\mu\mid m$, and the only information it uses on $m$ is that every prime dividing $F_1(p)$ exactly once divides $m$; hence $\mu$ must be squarefree. Moreover $F_1(p)=p(2p-1)(3p-2)$ is odd for odd $p$, so $\mu$ must be odd. Now $455=5\cdot7\cdot13$ is the least odd squarefree element of $\mathcal{S}$: by \Cref{lem:below-min}, any element of $\mathcal{S}$ is a multiple of an element of $\operatorname{Min}(\mathcal{S})$, and by \Cref{ex:min-S} the only odd element of $\operatorname{Min}(\mathcal{S})$ below $455$ is $315=3^2\cdot5\cdot7$, which is not squarefree. The residue $228$ is the inverse of $2$ modulo $455$, so $p\equiv228\pmod{455}$ means $455\mid 2p-1$, and the modulus $455^2$ ensures that $5$, $7$ and $13$ divide $2p-1$ only once. Distributing $5,7,13$ differently between $2p-1$ and $3p-2$ gives seven further residue classes modulo $455$, but the class of $228$ contains the least prime that works, namely $p=683$. Note that $683\equiv228\pmod{455}$ while $683\not\equiv228\pmod{455^2}$, so the congruence modulo $455^2$ is sufficient but not necessary.
\end{remark}

Requiring both $n$ and $n^2$ to be $\Delta$-primitive is much more restrictive. Below $20000$ there are only seven such $n$:
\[
84,\ 585,\ 1716,\ 7176,\ 7540,\ 16744,\ 17680 .
\]
In this range they are exactly the elements of $\operatorname{Min}(\mathcal{S})\cap\mathbb{N}(\Delta)$. The candidates $N_p$ do not help: $N_3=585$ is one of them, but $N_p\notin\mathbb{N}(\Delta)$ for $5\leq p\leq31$.

\begin{problem}
	\label{prob:double}
	Are there infinitely many $n$ such that both $n$ and $n^2$ are $\Delta$-primitive?
\end{problem}


	\section*{Acknowledgements}
	This work was partially supported by: Universidad Nacional de San Luis (UNSL), grants PI UNSL 03-1326, PROIPRO 03-2923, and PROINI 03-124; Consejo Nacional de Investigaciones Cient\'ificas y T\'ecnicas (CONICET), grant PIP 11220220100068CO; Gobierno de la Provicia de San Juan, Programa de Financiamiento ``Ìmpulsa Ciencia'' (Subprograma B).

\end{document}